\documentclass[a4paper,16pt]{article}
\usepackage{amsthm}
\usepackage[sumlimits]{amsmath}
\usepackage{amsfonts}
\usepackage{color}
\usepackage{mathrsfs}
\usepackage{colordvi,graphics,latexsym}
\allowdisplaybreaks

\def\MM#1{\boldsymbol{#1}}
\newtheorem{theorem}{Theorem}[section]

\DeclareMathOperator{\Id}{Id}

\newcommand{\p}{{\bf p}}
\newcommand{\q}{{\bf q}}

\newtheorem{corollary}[theorem]{Corollary}

\newcommand{\rem}[1]{}
\sffamily
\begin{document}

\title{Semi-geostrophic particle motion and 
exponentially accurate normal forms}
\author{Colin J Cotter\thanks{ Department of Mathematics,
Imperial College London, London SW7 2AZ, England, (e-mail: \texttt{
colin.cotter@imperial.ac.uk})}
\and 
Sebastian Reich\thanks{Institut f\"ur
Mathematik, Universit\"at Potsdam, Postfach 60 15 53, D-14415
Potsdam, Germany (e-mail:\texttt{sreich@math.uni-potsdam.de})}}

\maketitle

\begin{abstract}
We give an exponentially-accurate normal form for a Lagrangian
particle moving in a rotating shallow-water system in the
semi-geostrophic limit, which describes the motion in the region of an
exponentially-accurate slow manifold (a region of phase space for
which dynamics on the fast scale are exponentially small in the Rossby
number). We show how this result is
related to the variational asymptotics approach of \cite{oliver05};
the difference being that on the Hamiltonian side it is possible to
obtain strong bounds on the growth of fast motion away from (but near
to) the slow manifold. 

Our normal form approach extends to numerical approximations
\emph{via} backward error analysis, and extends to particle
methods for the shallow-water equations, where the
result shows that particles stay close to balance over long times
in the semi-geostrophic limit. 
\end{abstract}


%
%
%
\section{Introduction}
Semi-geostrophic approximations are obtained from the 
rotating shallow-water equations in the low Rossby
and Burger number limit. In \cite{oliver05}, a
systematic approach for deriving these models is set out, based on
looking for a near-identity change of coordinates for the solution
domain which makes the Lagrangian affine (linear in velocity),
resulting in a equation which only has slow dynamics.  
At each order in the Rossby number, there is some choice of
parameters available, and this choice can be crucial in obtaining a
well-posed PDE. After the change of coordinates, the equations of
motion can be obtained by applying Hamilton's variational principle to
the approximated Lagrangian.

In this paper we use exponentially-accurate normal form theory to
investigate the rotating shallow-water equations in the limit of low
Rossby and low Burger number. Exponentially-accurate normal
form theory was first popularised through the work of \cite{nek77:ham}
(very clearly reviewed and reformulated in
\cite{poes99:est_ell_equil}) on perturbed integrable Hamiltonian
ODEs. Also important in the literature, \cite{neis81:est_kol} gave an
exponentially-accurate estimate for systems with rapidly rotating
phase, and \cite{ben87:realisation} developed estimates for highly
oscillatory mechanical systems.  A first application of
exponentially-accurate normal form transformations in the context of
geophysical fluid dynamics can be found in \cite{wirshep00}. The
principle analytic obstruction to a result in our case is the
requirement that the free surface height remains sufficiently smooth;
hence, we restrict our result to individual particle trajectories
moving in a fluid with a specified free surface height (a problem
which has applications in tracking advected tracer particles: see the
summary for details), to toy ODE systems and to numerical solutions
obtained with the Hamiltonian particle-mesh method for the rotating
shallow-water equations.  The central tool of the normal form theory
is the choice of a suitable symplectic transformation for the particle
positions and momenta (this transformation may be more general than
the cotangent-lifted coordinate transformations used in
\cite{oliver05}), resulting in a Hamiltonian system where the small
coupling between fast and slow time-scales is more explicit. More
explicitly, if the model equations are solved with initial conditions
which are close to balance, then they will stay there for very long
time intervals. We also note that our results are more general than
the type provided in \cite{wir04}, where the existence of a slow
manifold up to exponentially small terms is shown.

The rest of the paper is organised as follows. In section
\ref{modelsystem} we introduce the model system of a particle moving
in a rapidly-rotating shallow-water flow, and discuss the ideas behind
geostrophy and semi-geostrophic approximations. In section \ref{Ham
form} we introduce the variational structure of these equations from
both of the Lagrangian and Hamiltonian viewpoints, and in section
\ref{higher} we discuss higher-order balanced models in the context of
the variational structure. This background then allows us to introduce
the exponentially-accurate normal form theory in section
\ref{exponential}, and we show how this theory extends to numerical
solutions obtained with symplectic time-stepping 
methods in section \ref{backward}.
Section \ref{particles} applies this theory to particle methods for
the rotating shallow-water equations. Section \ref{geostroph}
discusses geostrophic balance for these methods in the context of
exponentially-accurate normal form theory, whilst section
\ref{numerics} illustrates the theory with some numerical
results. Finally there is a summary and outlook in section
\ref{summary}.

%
\section{A model system}
\label{modelsystem}

The shallow-water equations (SWEs) on an $f$-plane are
\begin{eqnarray}
\frac{Du}{Dt} & = & +fv-g\mu _{x},\label{eq:u}\\
\frac{Dv}{Dt} & = & -fu-g\mu _{y},\label{eq:v}\\
\frac{D\mu }{Dt} & = & -\mu \left(u_{x}+v_{y}\right).\label{eq:mu}
\end{eqnarray}
Here $\mu =\mu \left(t,x,y \right)$ is the fluid
depth, $g$ is the gravitational constant, $f$
is twice the (constant) angular velocity of the reference plane, 
\begin{equation}
\frac{D}{Dt}\left(.\right)=\left(.\right)_{t}+u\left(.\right)_{x}
+v\left(.\right)_{y},\label{eq:Lag_def}
\end{equation}
is the Lagrangian or material time derivative, and subscripts denote
partial differentiation with respect to that variable. See \cite{salmon99}
for a derivation of the SWEs and their relevance to geophysical fluid
dynamics.

In non-dimensionalised variables, the SWEs can be written in the form
\begin{eqnarray}
{\rm Ro} \frac{Du}{Dt} & = & +v- \frac{{\rm B}}{\rm Ro} 
\mu _{x},\label{eq:nd_u}\\
{\rm Ro} \frac{Dv}{Dt} & = & -u- \frac{{\rm B}}{\rm Ro} 
\mu _{y},\label{eq:nd_v}\\
\frac{D\mu }{Dt} & = & -\mu \left(u_{x}+v_{y}\right),\label{eq:nd_mu}
\end{eqnarray}
where ${\rm Ro} = U/fL$ is the Rossby number and 
${\rm B} = (L_R/L)^2$ is the Burger number. Furthermore, 
$L_R = \sqrt{gH}/f$ is the Rossby radius of deformation,
$U$ is a typical advection velocity, $L$ is a typical length scale,
and $H$ is the mean fluid depth
for the problem under consideration. 
Semi-geostrophic theory is concerned with the limit of 
(\ref{eq:nd_u})-(\ref{eq:nd_mu}) with ${\rm Ro} = {\cal O}(\varepsilon)$, 
${\rm B} = {\cal O}(\varepsilon)$ as $\varepsilon \to 0$. For simplicity,
we use ${\rm Ro} = {\rm B} = \varepsilon$ from now on. 
See \cite{salmon99} for a derivation of non-dimensionalised equations 
and the semi-geostrophic scaling limit. 

Instead of investigating the full SWEs, we start in this paper
with the following simpler model problem. We assume that the layer
depth $\mu(\varepsilon t,x,y)$ is a given function of time $t$ and space 
${\bf x} = (x,y)^T$. We may then investigate the motion of a 
single fluid parcel with coordinates ${\bf q} = (q_x,q_y)^T$. 
The corresponding Newtonian equations of motion are 
given by 
\begin{eqnarray}
\varepsilon \frac{d^2 q_x}{dt^2} &=& + \frac{d q_y}{dt} - 
\mu_{q_x}(\varepsilon t,q_x,q_y), \label{eq:ddot_x}\\
\varepsilon \frac{d^2 q_y}{dt^2} &=& - \frac{d q_x}{dt} - 
\mu_{q_y}(\varepsilon t,q_x,q_y), \label{eq:ddot_y}
\end{eqnarray}
We next rescale time and introduce the new time-scale $\tau = 
\varepsilon \,t$. We denote time derivatives with respect
to $\tau$ by overdot, \emph{e.g.}, $dq_x/d\tau = \dot{q}_x$. We also
write the second order equations (\ref{eq:ddot_x})-(\ref{eq:ddot_y}) 
as a system of first order equations by introducing the momentum 
${\bf p} = (p_x,p_y)^T$, \emph{i.e.},
\begin{eqnarray}
\dot{p}_x &=& + p_y - \varepsilon \,\mu_{q_x}(\tau,q_x,q_y), 
\label{eq:dot_u}\\
\dot{p}_y &=& - p_x - \varepsilon \,\mu_{q_y}(\tau,q_x,q_y), 
\label{eq:dot_v}\\
\dot{q}_x &=& p_x, \label{eq:dot_x}\\
\dot{q}_y &=& p_y. \label{eq:dot_y}
\end{eqnarray}
These equations may be expressed in a more compact form:
\begin{eqnarray}
\dot{\bf p} &=& J_2 {\bf p} - \varepsilon \,\nabla_{\bf q} 
\mu(\tau,{\bf q}), 
\qquad  J_2 = \begin{pmatrix}0 & 1 \\-1 & 0 \\
\end{pmatrix}, \label{eq:dot_p}\\ 
\dot{\bf q} &=& {\bf p}, \label{eq:dot_q}
\end{eqnarray}
which is the formulation we will work with in this paper. A
discussion of the relationship between the Hamiltonian structure of this
system and the Hamiltonian structure for the full SWEs
is given in \cite{bokhove05:_parcel_euler_lagran}.

The intuitive idea behind the semi-geostrophic approximation
\cite{salmon99} is that the solutions of
(\ref{eq:dot_p})-(\ref{eq:dot_q}) consist of inertial oscillations
with period $T_I = 2\pi = {\cal O}(\varepsilon^0)$ and slow
geostrophically balanced motion on a (slow) time-scale $T_G = {\cal
O}(\varepsilon)$. Furthermore, intuition suggests that the inertial
oscillations are primarily governed by the linear equations
\begin{eqnarray}
\dot{\bf p} &=& J_2 {\bf p}, \label{eq:dot_p_in}\\ 
\dot{\bf q} &=& {\bf p}, \label{eq:dot_q_in}
\end{eqnarray}
while the slow, geostrophically balanced, parcel dynamics is
characterised by the reduced (nonlinear) system
\begin{eqnarray}
{\bf 0} &=& J_2 {\bf p} - \varepsilon \,\nabla_{\bf q} 
\mu(\tau,{\bf q}), \label{eq:dot_p_gb}\\ 
\dot{\bf q} &=& {\bf p}, \label{eq:dot_q_gb}
\end{eqnarray}
or, equivalently,
\begin{equation} \label{eq:gsa}
\dot{\bf q} = -\varepsilon \,J_2 \nabla_{\bf q} \mu(\tau,{\bf q}).
\end{equation}

%
\subsection{Hamiltonian and variational formulations}
\label{Ham form}

We start with a Hamiltonian formulation of the system 
(\ref{eq:dot_p})-(\ref{eq:dot_q}). As a further simplification
we assume that the layer-depth $\mu$ is time-independent and
introduce the potential energy function $V({\bf q}) := 
\mu({\bf q})$ to make the link to classical mechanics more
transparent.

We rewrite the equations (\ref{eq:dot_p})-(\ref{eq:dot_q}) 
using a non-canonical symplectic structure
operator
\begin{equation} \label{eq:J4}
J = \begin{pmatrix}
J_2 & -I_2 \\
I_2 & 0_2 \\
\end{pmatrix} \in \mathbb{R}^{4\times 4},
\end{equation}
so that
\begin{equation} \label{eq:abstract_ham}
\dot{\bf z} 
= J \nabla_{\bf z} H_0({\bf z}),
\end{equation}
with Hamiltonian
\begin{equation} \label{eq:hamiltonian_0}
H_0({\bf z}) = K({\bf p}) + \varepsilon V({\bf q}), 
\quad K({\bf p}) = \frac{1}{2}\p^T\p,
\end{equation}
and phase space variable ${\bf z} = ({\bf p}^T,{\bf q}^T)^T \in
\mathbb{R}^4$.

Another approach is to rewrite the system of first-order equations
as a second-order equation
\begin{equation} \label{eq:2nd}
\ddot{\bf q} - J_2 \dot{\bf q} + \varepsilon \nabla_{\bf q} 
V({\bf q}) = {\bf 0},
\end{equation}
and to note that (\ref{eq:2nd}) is the Euler-Lagrange equation
for the Lagrangian functional
\begin{equation} \label{eq:lvp1}
{\cal L} = \int d\tau \left[ \frac{1}{2} \|\dot{\bf q}\|^2 +
\frac{1}{2} {\bf q}^T J_2 \dot{\bf q} - \varepsilon V({\bf q})\right].
\end{equation}

The reduced model (\ref{eq:gsa}) with $\mu = V$ is also Hamiltonian
with phase space ${\bf q} \in \mathbb{R}^2$, structure matrix $J_2^T$, and 
Hamiltonian function $H_g({\bf q}) = \varepsilon V({\bf q})$. The associated
Lagrangian functional is given by
\begin{equation} \label{eq:lvp2}
{\cal L}_{\rm g} = \int d\tau \left[ 
\frac{1}{2} {\bf q}^T J_2 \dot{\bf q} - \varepsilon V({\bf q})\right].
\end{equation}
Note that (\ref{eq:lvp2}) differs from (\ref{eq:lvp1}) by the
missing kinetic energy term $\|\dot{\bf q}\|^2/2$. 

%
\subsection{Higher order balance and the `semi-geostrophic' 
approximation}
\label{higher}

Clearly, one would like to derive more sophisticated reduced
models; a systematic approach has recently been given in
\cite{oliver05}. In this section we summarize the main results
from \cite{oliver05} before we develop our novel approach in 
section \ref{exponential}.

All known derivations of balanced (semi-geostrophic) models start from
the assumption that
\begin{equation} \label{eq:boundedderivative}
\dot{\bf q} = -\varepsilon J_2 \nabla_{\bf q} V({\bf q})
+{\cal O}(\varepsilon^2) = {\cal O}(\varepsilon).
\end{equation}
Under this assumption, we may formally collect terms of equal order
in $\varepsilon$ and rewrite (\ref{eq:lvp1}) as 
\begin{equation} \label{eq:lvpexpan}
{\cal L} = \varepsilon \int d\tau \left[ L_0 + \varepsilon L_1 \right],
\quad \mbox{where} 
\quad L_0 = \frac{1}{2\varepsilon} {\bf q}^T J_2 \dot{\bf q} -
V({\bf q}), \quad L_1 = \frac{1}{2\varepsilon^2} \|\dot{\bf q}\|^2.
\end{equation}
We now introduce a coordinate transformation $\psi_\varepsilon:
\mathbb{R}^2 \to \mathbb{R}^2$ and transformed coordinates
${\bf q}_\varepsilon$ via
\begin{equation} \label{eq:q-trans}
{\bf q} = \psi_\varepsilon ({\bf q}_\varepsilon) = 
{\bf q}_\varepsilon + \varepsilon {\bf F}_1({\bf q}_\varepsilon) 
+ {\cal O}(\varepsilon^2).
\end{equation}
Following \cite{oliver05} (but note our different notation), we set
\begin{equation} \label{eq:q-geostrophic}
{\bf F}_1({\bf q}_\varepsilon) = 
-\frac{1}{2}\nabla_{\bf q} V({\bf q}_\varepsilon) + \lambda
\nabla_{\bf q}V({\bf q}_\varepsilon) = -\frac{1}{2\varepsilon} J_2
\dot{\bf q}_\varepsilon 
+ \lambda \nabla_{\bf q} V({\bf q}_\varepsilon) + {\cal O}(\varepsilon).
\end{equation}
Hence we obtain a Lagrangian function
\begin{equation} \label{eq:trans-lagr}
{\cal L}_1 = \varepsilon \int d\tau \left[ \bar L_0 + \varepsilon
\bar L_1 + {\cal O}(\varepsilon^2) \right]
\end{equation}
in the transformed variable ${\bf q}_\varepsilon$, where
\begin{equation}
\bar L_0 = \frac{1}{2\varepsilon} {\bf q}_\varepsilon^T J_2 
\dot{\bf q}_\varepsilon - V({\bf q}_\varepsilon), \qquad
\bar L_1 = \varepsilon^{-1} \left[ \frac{1}{2}+\lambda \right] 
\nabla_{\bf q} V({\bf q}_\varepsilon)^T J_2 \dot{\bf q}_\varepsilon
- \lambda \|\nabla_{\bf q} V({\bf q}_\varepsilon)\|^2.
\end{equation} 
Upon dropping ${\cal O}(\varepsilon^3)$ terms in (\ref{eq:trans-lagr}),
we finally obtain the transformed Lagrangian
\begin{equation} 
{\cal L}_1 = \int d\tau \left[ 
\frac{1}{2} {\bf q}_\varepsilon^T J_2 
\dot{\bf q}_\varepsilon - \varepsilon V({\bf q}_\varepsilon)
+  \varepsilon \left[ \frac{1}{2}+\lambda \right] 
\nabla_{\bf q} V({\bf q}_\varepsilon)^T J_2 \dot{\bf q}_\varepsilon
- \varepsilon^2 \lambda \|\nabla_{\bf q} V({\bf q}_\varepsilon)\|^2
\right].
\end{equation}

In the context of this paper, the choice $\lambda = -1/2$ is of particular
interest. Note that this choice corresponds to Salmon's 
large-scale semi-geostrophic approximation (see \cite{salmon99,oliver05}). 
The associated reduced equations of motion can be derived 
from the Lagrangian functional 
\begin{equation} \label{eq:lvp3}
{\cal L}_{\rm lsg}  = 
 \int d\tau \left[ \frac{\varepsilon^2}{2} \|\nabla_{\bf q} 
V({\bf q}_\varepsilon)
\|^2 + \frac{1}{2} {\bf q}^T_\varepsilon J_2 \dot{\bf q}_\varepsilon 
- \varepsilon  V({\bf q}_\varepsilon)\right].
\end{equation}
The associated Euler-Lagrange equations are explicitly given by
\begin{equation} \label{eq:igsa}
J_2 \dot{\bf q}_\varepsilon - 
\varepsilon \nabla_{\bf q} V({\bf q}_\varepsilon) + 
\frac{\varepsilon^2}{2} \nabla_{\bf q} \| \nabla_{\bf q} 
V({\bf q}_\varepsilon)\|^2 = 
{\bf 0}.
\end{equation}
These equations are canonical with Hamiltonian
\begin{equation}
H_{\rm lsg}({\bf q}_\varepsilon) 
= \varepsilon V({\bf q}_\varepsilon) - \frac{\varepsilon^2}{2}
\|\nabla_{\bf q} V({\bf q}_\varepsilon)\|^2,
\end{equation}
and structure matrix $J_{\rm lsg} = J_2^T$.

A number of questions arises at this stage.

\begin{enumerate}

\item The assumption (\ref{eq:boundedderivative}) is essential 
for the derivation of the reduced Lagrangian (\ref{eq:trans-lagr}).
However, even if (\ref{eq:boundedderivative}) holds at the initial 
time, how can one ensure that it holds true for solutions of the full
equations (\ref{eq:2nd}) over finite (or infinite) time intervals?

\item Following \cite{oliver05}, one can formally generalize the
coordinate transformation to obtain higher-order balanced equations.
Can this process be continued to arbitrary precision and what 
rigorous error bounds in terms of $\varepsilon$ can be obtained?

\item What role is played by the parameter $\lambda$?

\end{enumerate}

We will provide an answer to these questions from the perspective
of Hamiltonian perturbation theory in the following sections.
We will first restrict the discussion to the simple model system
(\ref{eq:dot_p})-(\ref{eq:dot_q}) before considering a generalization to
particle methods for the SWEs in section \ref{particles}.

We also wish to point the reader to the reports 
\cite{gottwald1,gottwald2}, which attempt to answer the above
questions based on higher-order coordinate transformations
\begin{equation} \label{eq:marcel1}
{\bf q}_\varepsilon = {\bf q} + \varepsilon {\bf F}_1 + \varepsilon^2 
{\bf F}_2 + \cdots + \varepsilon^{m} {\bf F}_m + 
{\cal O}(\varepsilon^{m+1})
\end{equation}
and associated higher-order expansions 
\begin{equation} \label{eq:marcel2} 
{\cal L}_m = \varepsilon \int d\tau \left[ \bar L_0 + \varepsilon
\bar L_1 + \cdots + \varepsilon^m 
\bar L_m + {\cal O}(\varepsilon^{m+1}) \right]
\end{equation}
of the Lagrangian functional for $m\ge 1$.

%
\subsection{Exponentially accurate normal forms}
\label{exponential}

In this section we view the finite dimensional system 
(\ref{eq:abstract_ham}) from the
Hamiltonian side and look for symplectic coordinate transformations to
a normal form. 


The precise aim of the normal form transformation is to find a symplectic
(with respect to the structure matrix $J$) change of
coordinates, which transforms the system to a form in which the fast
and slow variables are completely separated up to a small remainder term. 
As the transformation is symplectic, the equations of motion can be 
obtained simply by substituting the change of coordinates into 
the Hamiltonian and writing down Hamilton's equations. This normal
form strategy is well-developed and we only summarise the main
steps. See, \emph{e.g.}, \cite{Neishtadt84} for the first derivation of
an exponential estimate and \cite{cotter04} for an application
to systems of type (\ref{eq:abstract_ham}).

Our aim is to find a near-identity change of coordinates $\Psi_n$ so
that
\begin{equation}
H_n = H_0\circ\Psi_n = K + \varepsilon G_n + \varepsilon^{n+1}R_n,
\end{equation}
where
\begin{equation}
\{G_n,K\} = 0,
\end{equation}
with $\{\cdot,\cdot\}$ being the Poisson bracket obtained from
$J$.

We define $\Psi_{n+1}$ recursively by writing 
\begin{equation}
\Psi_{n+1} = \Psi_{n}\circ\Phi_{1,\varepsilon^{n+1} F_{n+1}}, \quad
\Psi_0 = \Id,
\end{equation}
where $\Phi_{1,\varepsilon^{n+1} F_{n+1}}$ is the time-1 flow map of 
the Hamiltonian vector field produced by the Hamiltonian function 
$H = \varepsilon^{n+1} F_{n+1}$. Recursive substitution gives
\begin{equation}
H_{n+1} = H_{n}\circ\Phi_{1,\varepsilon^{n+1} F_{n+1}}
= K + \varepsilon G_{n} + \varepsilon^{n+1} (R_n + \{K,F_{n+1}\})
+ \mathcal{O}(\varepsilon^{n+2}),
\end{equation}
and so we need to choose $F_{n+1}$ so that
\begin{equation}
\label{homol eqn}
R_n + \{K,F_{n+1}\} = \bar{R}_n,
\end{equation}
where 
\begin{equation}
\{\bar{R}_n,K\} = 0.
\end{equation}
Equation (\ref{homol eqn}) is known as the \emph{homological
equation}. The solution $F_{n+1}$, and the \emph{non-resonant part}
$\bar{R}_n$ of $R_n$ can be calculated \emph{via}
\begin{eqnarray}
\bar{R}_n & = & \frac{1}{T}\int_0^T d\tau 
\, R_n\circ\Phi_{\tau,K}, \label{R n+1} \\
F_{n+1} & = & \frac{1}{T}\int_0^T d\tau \,\tau 
(R_n-\bar{R}_n)
\circ\Phi_{\tau,K}, \label{F n+1}
\end{eqnarray}
where $\Phi_{\tau,K}$ is the time-$\tau$ flow-map of the Hamiltonian vector
field produced by $K$, which is time-periodic with period $T=2\pi$.
To close the recursion we define
\begin{equation}
G_{n+1} = G_n + \varepsilon^n \bar{R}_n
\end{equation}
and terms of size $\mathcal{O}(\varepsilon^{n+2})$ and smaller 
are collected in a new residual function $R_{n+1}$. 

When the potential $V$ is analytic and bounded on an appropriate 
compact set in the complex plane, we can choose to make the number
of iterations scale with $1/\varepsilon$ and obtain an
exponentially-accurate normal form as summarised in the following 
theorem:
\begin{theorem}
\label{normal form theorem}
Let $\mathscr{K}$ be some compact subset of phase space containing the
solution trajectory, and let $V$ be analytic in $B_r\mathscr{K}$ ( the
union of complex balls of radius $r$ with centres in $\mathscr{K}$).
Then there exists a near-identity change of coordinates
$\Psi_{\varepsilon}$ such that
\begin{equation}
\tilde{H} = 
H_0\circ\Psi_{\varepsilon} = K+\varepsilon \,G
+e^{-c/\varepsilon}\,R,
\end{equation}
where
\begin{equation}
\{G,K\}=0,
\end{equation}
$G$ and $R$ have $\epsilon$-independent bounds on $\mathscr{K}$ and
$c>0$ is some constant.
\end{theorem}
\begin{proof}
Essentially the method is to obtain a Cauchy estimate
for each remainder $R_n$, and to truncate the above algorithm after
$n$ steps where $n$ scales with $1/\varepsilon$. Defining the
norm on functions
\begin{equation}
|f|_r = \sup_{\MM{z}\in B_r\mathscr{K}}|f(\MM{z})|,
\end{equation}
we obtain iterative estimates for $\bar{R}_n$ and $F_{n+1}$ from
equations (\ref{R n+1}-\ref{F n+1}) \cite{poes99:est_ell_equil}:
\begin{equation}
|\bar{R}_n|_r \leq |R_n|_r, \qquad 
|F_{n+1}|_r \leq 2\pi|R_n|_r.
\end{equation}
Simple adding and subtracting of terms gives
\begin{eqnarray}
  \nonumber H_n\circ\Phi_{1,\epsilon^{n+1}F_{n+1}} &=& K
  + [\epsilon G_n+\epsilon^{n+1}\{K,F_{n+1}\}+\epsilon^{n+1}R_n] \\
  \nonumber & & + K\circ\Phi_{1,\epsilon^{n+1}F_{n+1}} -K-\epsilon^{n+1}\{K,F_{n+1}\} \\
  \nonumber & & + \epsilon(G_n\circ\Phi_{1,\epsilon^{n+1}F_{n+1}}-G_n) \\
  & & + \epsilon^{n+1}(R_n\circ\Phi_{1,\epsilon^{n+1}F_{n+1}}-R_n). \label{H_{n+1}}
\end{eqnarray}
The term in the square brackets is $\epsilon G_{n+1}$, and it 
remains to estimate this term, and the remainder. If we choose
$\delta$ so that $r-(n+1)\delta>0$ and $H_n$ is analytic in 
$B_{r-n\delta}\mathcal{K}$ then
\begin{eqnarray}
\nonumber |G_{n+1}|_{r-(n+1)\delta} &=& 
|G_n+\epsilon^n\overline{R}_n|_{r-(n+1)\delta}, \\
&\leq& |G_n|_{r-n\delta}+\epsilon^n|R_n|_{r-n\delta}, \label{G n estimate}
\end{eqnarray}
having used the estimate for $\overline{R}_n$. Estimates for the
remainder $R_{n+1}$ can be obtained using the mean-value theorem
combined with a Cauchy estimate for the gradient. To estimate the
second line in equation (\ref{H_{n+1}}) note that
\begin{eqnarray} \nonumber
 & &
\mathcal{K}\circ\Phi_{1,\epsilon^{n+1}F_{n+1}}-
\mathcal{K}-\epsilon^{n+1}\{\mathcal{K},g_{n+1}\} \\ \nonumber
& = & \epsilon^{n+1}\int_0^1\{\mathcal{K},F_{n+1}\}\circ
\Phi_{t,\epsilon^{n+1}F_{n+1}}dt -\epsilon^{n+1}\{\mathcal{K},F_{n+1}\}, \\
\nonumber
& = & \epsilon^{n+1}\int_0^1\{\mathcal{K},F_{n+1}\}\circ\Phi_{t,\epsilon^{n+1}F_{n+1}}-\{\mathcal{K},F_{n+1}\}dt, \\ 
& = & \epsilon^{n+1}\int_0^1(\overline{R}_n-R_n)\circ\Phi_{t,\epsilon^{n+1}F_{n+1}}-(\overline{R}_n-R_n)dt,
\end{eqnarray}
which may be bounded using the mean-value theorem and a Cauchy estimate
by
\begin{eqnarray}
  \nonumber & & 
  \nonumber |\mathcal{K}\circ\Phi_{1,\epsilon^{n+1}F_{n+1}}-\mathcal{K}-\epsilon^{n+1}\{\mathcal{K},g_{n+1}\}|_{r-(n+1)\delta} \\
  & \leq & 
  \nonumber \epsilon^{n+1}\int_0^1\left|(\overline{R}_n-R_n)\circ\Phi_{t,\epsilon^{n+1}F_{n+1}}-(\overline{R}_n-R_n)\right|_{r-(n+1)\delta}dt, \\
  \nonumber & \leq & \epsilon^{2(n+1)}|\nabla(\overline{R}_n-R_n)|_{r-(n+1)\delta}
  |F_{n+1}|_{r-(n+1)\delta} \\
  & \leq & \frac{2\pi\epsilon^{2(n+1)}}{\delta}|\overline{R}_n-R_n|_{r-n\delta}
  |R_n|_{r-(n+1)\delta}, \nonumber \\
  & \leq & \frac{2\pi\epsilon^{2(n+1)}}{\delta}|R_n|_{r-n\delta}
  |R_n|_{r-n\delta}. \label{R est 1}
\end{eqnarray}
Similar estimates give
\begin{eqnarray}
  \label{R est 2} \epsilon |G_n\circ\Phi_{1,\epsilon^{n+1}F_{n+1}}-G_n|_{r-(n+1)\delta} & \leq &
  \frac{2\pi\epsilon^{n+2}}{\delta}
  |G_n|_{r-n\delta}|R_n|_{r-n\delta}, \\
  \epsilon^{n+1}|R_n\circ\Phi_{1,\epsilon^{n+1}F_{n+1}}-R_n|_{r-(n+1)\delta} & \leq &
  \label{R est 3} \frac{2\pi\epsilon^{2(n+1)}}{\delta}|R_n|_{r-n\delta}|R_n|_{r-n\delta}.
\end{eqnarray}
Suppose we have:
\begin{equation}
|G_n|_{r-n\delta}\leq c_0\sum_{i=1}^n\left(\frac{c_1\epsilon}{\delta}\right)^i,
\quad |R_n|_{r-n\delta}\leq c_0\left(\frac{c_1}{\delta}\right)^n,
\end{equation}
where $c_0$ is a positive constant, and $c_1 = 4\pi c_0$,
then we get
\begin{equation}
|G_{n+1}|_{r-(n+1)\delta}\leq
c_0\sum_{i=1}^{n+1}\left(\frac{c_1\epsilon}{\delta}\right)^i,
\end{equation}
from equation (\ref{G n estimate}), and
\begin{eqnarray} \nonumber
|R_{n+1}|_{r-(n+1)\delta} & = & \frac{2\pi\epsilon^{n+1}}{\delta}|R_n|_{r-n\delta}
  |R_n|_{r-n\delta} +   \frac{2\pi}{\delta}
  |G_n|_{r-n\delta}|R_n|_{r-n\delta} \\ \nonumber
 & & \qquad + \frac{2\pi\epsilon^{n+1}}{\delta}|R_n|_{r-n\delta}
|R_n|_{r-n\delta}, \\ 
& \leq & c_0\left(\frac{c_1}{\delta}\right)^n\left(
\frac{4\pi\epsilon^{n+1}}{\delta}c_0\left(\frac{c_1}{\delta}\right)^n
  +   \frac{2\pi}{\delta}c_0\sum_{i=1}^n\left(
\frac{c_1\epsilon}{\delta}\right)^i\right) , 
\end{eqnarray}
from equations (\ref{R est 1}-\ref{R est 3}). If
\begin{equation}
\label{half cond}
\frac{c_1\epsilon}{\delta}<\frac{1}{2},
\end{equation}
then $|G_{n+1}|_{r-(n_1)\delta}\leq c_0$, and
\begin{equation}
|R_{n+1}|_{r-(n+1)\delta} \leq
c_0\left(\frac{c_1}{\delta}\right)^n\frac{2\pi c_0}{\delta}
\left(2\epsilon\frac{1}{2^n} + 1
\right) \leq c_0\left(\frac{c_1}{\delta}\right)^{n+1}.
\end{equation}
These estimates are true for all $n$ by induction if we set $c_0 =
|V|_r$. Finally we choose $n$ to scale with $\epsilon$ from
the formula
\begin{equation}
n = \left[
\frac{2r}{ec_1\epsilon}
\right]
\end{equation}
where the square brackets indicate rounding up to the nearest integer,
and set $\delta=r/2n$. Equation (\ref{half cond}) can then be satisfied, 
and we get
\begin{equation}
\epsilon^{n+1}|R_n|_{r/2} \leq \epsilon c_0 (\frac{\epsilon c_1}{\delta})^n
\leq \epsilon c_0\left(
\frac{1}{e}\right)^{2/e c_1\epsilon}
= \epsilon c_0e^{-c\epsilon}
\end{equation}
where $c = e c_1/2$.
\end{proof}

Let us denote the transformed variables by ${\bf z}_\varepsilon = 
(\p^T_\varepsilon,\q^T_\varepsilon)^T$. The transformed variables
are defined by $\Psi_\varepsilon({\bf z}_\varepsilon) = {\bf z}$. 
We can then obtain an estimate 
for $K(\p_\varepsilon(\tau))$ in these coordinates by
\begin{eqnarray} \nonumber
|K(\p_\varepsilon(\tau))-K(\p_\varepsilon(0))| & < & |\tau| \left|
\frac{d}{d\tau}K(\p_\varepsilon(\tau)) \right|, \\ \nonumber
 & < & |\tau||\{K,\tilde{H}\}|, \\ \nonumber
& < & |\tau|e^{-c/\varepsilon}|\{K,R\}|, \\ 
& < & d e^{-c/2\varepsilon}, \label{eq:KE_drift}
\end{eqnarray}
for $|\tau|<e^{c/2\varepsilon}$ and some constant $d>0$. This means that
$\|\p_\varepsilon(\tau)\|^2$ stays almost constant for very long 
time intervals.

\begin{corollary} \label{corollary1}
Let us assume that the momenta ${\bf p} = \dot{\bf {q}}$
satisfy
\begin{equation} \label{eq:initialb}
{\bf p}(0) = -\varepsilon J_2 \nabla_{\bf q}V({\bf q}(0)) + 
{\cal O}(\varepsilon^2)
\end{equation}
at initial time $\tau = 0$, then 
\begin{equation}
{\bf p}(\tau) = -\varepsilon J_2 \nabla_{\bf q}V({\bf q}(\tau)) + 
{\cal O}(\varepsilon^2)
\end{equation}
for all $|\tau|<e^{c/2\varepsilon}$.
\end{corollary}
\begin{proof} To first order in $\varepsilon$ the slow manifold is
  determined by ${\bf p} = -\varepsilon J_2 \nabla_{\bf q}V({\bf q})$
  and, in terms of the transformed coordinates, by ${\bf
    p}_\varepsilon = {\cal O}(\varepsilon^2)$. Hence
  (\ref{eq:initialb}) implies that $K({\bf p}_\varepsilon(0)) = {\cal
    O}(\varepsilon^4)$. Finally (\ref{eq:KE_drift}) yields the desired
  result.
\end{proof}
This means that provided the system is within
$\mathcal{O}(\varepsilon^2)$ of the geostrophically balanced state
initially, it will stay there for exponentially long time intervals.
It is \emph{not} necessary for the system to be exactly
geostrophically balanced for the the result to hold. 

If $\p_\varepsilon(0)=0$, (\emph{i.e.}, we start on the slow manifold)
then there are no rapid oscillations due to $\p_\varepsilon$, and we
get the symplectic slow equation given in the following corollary:
\begin{corollary}
If $\p_\varepsilon(0)=0$ then 
\begin{equation} \label{eq:gse}
J_2\dot{\q}_\varepsilon =  \varepsilon \nabla_{
\q} G ({\bf 0},{\bf q}_\varepsilon)+ \mathcal{O}(e^{-c/2\varepsilon}),
\end{equation}
for all $|\tau|<e^{c/2\varepsilon}$.
\end{corollary}
\begin{proof}
  If $\p_\varepsilon(0)=0$, then theorem \ref{normal form theorem}
  says that $\p_\varepsilon(t)$ stays exponentially small for
  exponentially long time intervals and so
  \begin{equation}
    \label{q hamil}
    \dot{\q}_\varepsilon = 
    \varepsilon \, \nabla_{\p} G({\bf p}_\varepsilon,
    {\bf q}_\varepsilon)|_{\p_\varepsilon={\bf 0}} +  
    \mathcal{O}(e^{-c/2\varepsilon}),
  \end{equation}
  for all $|\tau|<e^{c/2\varepsilon}$.
  We can obtain the slow equations in symplectic form by 
  examining the $\p_\varepsilon$ equation and ignoring exponentially small
  terms, \emph{i.e.}, 
\begin{equation}
{\bf 0} = \dot{\p}_\varepsilon = 
J_2 \nabla_{\p}G({\bf p}_\varepsilon,{\bf q}_\varepsilon)
|_{\p_\varepsilon={\bf 0}} 
- \nabla_{\q} G({\bf p}_\varepsilon,{\bf q}_\varepsilon)
|_{\p_\varepsilon ={\bf 0}} ,
\end{equation}
and then substituting equation (\ref{q hamil}) to make
\begin{equation}
J_2\dot{\q}_\varepsilon = \varepsilon \nabla_{
\q} G ({\bf p}_\varepsilon,{\bf q}_\varepsilon) |_{\p_\varepsilon=0} 
=  \varepsilon \nabla_{
\q} G ({\bf 0},{\bf q}_\varepsilon),
\end{equation}
up to terms exponentially small in $\varepsilon$.  
\end{proof}
We note that
(\ref{eq:gse}) is in canonical form, \emph{i.e.}, the equations are
Hamiltonian with structure matrix $J_{\rm lsg} = J_2^T$ and
Hamiltonian $\hat H_{\rm lsg}({\bf q}_\varepsilon) = \varepsilon\,
G({\bf 0},{\bf q}_\varepsilon)$. When the system is \emph{not}
initialised on the slow manifold then the dynamics is dominated by
this slow equation but also contains fast components (which remain
approximately of the same magnitude for exponentially long time
intervals).

In section \ref{modelnumerics} we will conduct a simple experiment
that will allow us to verify the above estimates numerically. The
particular set-up will be chosen such that ${\bf p} = {\bf
  p}_\varepsilon$ at the initial time $\tau = 0$ and final time $\tau
= T$. Then $|K({\bf p}(0)) - K({\bf p}(T))|$ should go to zero
exponentially fast as $\varepsilon \to 0$ (discounting numerical
round-off errors). We will show in section \ref{backward} that such
exponentially small terms can indeed be observed using a symplectic
time stepping method (see, {\it e.g.}, \cite{HLW02,LR2005} for a
discussion of symplectic time stepping methods for Hamiltonian ODEs).

In the appendix we calculate the slow equation (\ref{eq:gse}) 
to second-order by
applying the iterative algorithm and obtain
\begin{equation} \label{eq:igsa_ham}
J_2\dot{\q}_\varepsilon = 
\varepsilon\nabla_{\q} V(\q_\varepsilon) - 
\frac{\varepsilon^2}{2}\nabla_{\q} \|\nabla_{\q} 
V(\q_\varepsilon)\|^2 + \mathcal{O}(\varepsilon^3).
\end{equation}
Note that the significant terms in (\ref{eq:igsa_ham}) 
are equivalent to the `large-scale semi-geostrophic' equation (\ref{eq:igsa}).
In other words, (\ref{eq:gse}) provides a higher order generalisation
of the `large-scale semi-geostrophic' theory in the 
context of our simple model system. More precisely, Hamiltonian
normal form theory naturally leads to $\lambda = -1/2$.

%
\subsection{Extension to non-autonomous systems}
\label{time-dependent}
To recover results for equations (\ref{eq:ddot_x})-(\ref{eq:ddot_y}) we
need to extend the results from section \ref{exponential} 
to systems with time-dependent potentials where
the Hamiltonian takes the form
\begin{equation}
H_0({\bf z},\tau) = K({\bf p}) + \varepsilon V({\bf q},\tau).
\end{equation}
We apply the 
standard technique of extending the phase space by two extra 
conjugate variables $s$ and $e$, and write a new Hamiltonian
\begin{equation}
\bar{H}_0({\bf z},e,s,\tau) = K(\p) + \varepsilon V(\q,s) - e,
\end{equation}
so that the equations of motion for $e$ and $s$ are
\begin{equation}
\dot{s} = 1, \qquad \dot{e} = \varepsilon \frac{\partial V}{\partial s}
(\q,s),
\end{equation}
and we recover the original equations. The new Hamiltonian $\bar{H}_0$
is now autonomous and in a suitable form to apply theorem \ref{normal
form theorem} with minor modifications (to account for the two
extra variables), \emph{i.e.}, there exists a near-identity symplectic
change of coordinates
$(\p,\q,e,s)\to({\p}_\varepsilon,{\q}_\varepsilon,{e}_\varepsilon,
{s}_\varepsilon)$ such that
$K({\p}_\varepsilon)$ is nearly preserved over exponentially long time
intervals.

We will make use of this extension to non-autonomous systems again
in section \ref{numerics}.

%
\subsection{Numerical methods and backward error analysis}
\label{backward}

The Hamiltonian equations (\ref{eq:dot_p})-(\ref{eq:dot_q}) can be
solved numerically using a symplectic method. As an example, we give
the following splitting method. We rewrite the Hamiltonian
(\ref{eq:hamiltonian_0}) as a sum
\begin{equation}
H_0 = \frac{1}{2} K + \varepsilon V + 
\frac{1}{2} K
\end{equation}
and note that each entry has an exact flow map, which we denote
by $\Phi_{\tau,K/2}$ and $\Phi_{\tau,\varepsilon V}$, respectively.
A second-order symplectic method is now given by the composition
of these flow maps \cite{LR2005}, \emph{i.e.},
\begin{equation} \label{eq:time_stepping}
{\bf z}^{n+1} = M_{\Delta \tau}({\bf z}^n), \qquad
M_{\Delta \tau} = 
\Phi_{\Delta \tau/2,K} \circ \Phi_{\Delta \tau,\varepsilon V} \circ
\Phi_{\Delta \tau/2,K} .
\end{equation}

Using backward error analysis, it can be shown that this method
is equivalent to the exact solution of a modified Hamiltonian problem
up to terms exponentially small in $\Delta \tau$ \cite{BenGio94,Reich99}. 
More specifically, there exists a modified Hamiltonian
\begin{equation} \label{eq:mod_ham}
H_{\Delta \tau}({\bf z}) = H_0({\bf z}) 
+ \varepsilon \,\Delta \tau^2 P({\bf z},
\Delta \tau)
\end{equation}
such that
\begin{equation} \label{eq:modhamest}
\| M_{\Delta \tau}({\bf z}) - \Phi_{\Delta \tau,
H_{\Delta \tau}}({\bf z})
\| \le c_1\, e^{-c_2/\Delta \tau},
\end{equation}
where $c_1,c_2 >0$ are appropriate constants. See \cite{HLW02,LR2005}
for a survey of backward error analysis results for Hamiltonian systems. 

The modified Hamiltonian (\ref{eq:mod_ham}) allows us 
to discuss the preservation of geostrophic balance
under the time-stepping method (\ref{eq:time_stepping}). All one has
to do is to apply the normal form transformation from the 
previous section to the modified Hamiltonian (\ref{eq:mod_ham})
and to generalise theorem \ref{normal form theorem}. 
In particular, the following theorem implies the near conservation of 
geostrophic balance under the numerical method (\ref{eq:time_stepping})
over exponentially long times. 

\begin{theorem}
  \label{backward error analysis}
  Consider a symplectic integrator 
  \begin{equation} 
   {\bf z}^{n+1} = M_{\Delta t}({\bf z}^n), \qquad {\bf z}^n = 
   (({\bf q}^n)^T,({\bf p}^n)^T)^T,
  \end{equation}
  such as (\ref{eq:time_stepping}). 
  Then, for a sufficiently small time step
  $\Delta \tau$, there exists a symplectic change of coordinates 
  $\Psi_{\varepsilon,\Delta \tau}$ such
  that 
\begin{equation}
\| \Psi_{\varepsilon,\Delta \tau}^{-1} \circ 
M_{\Delta \tau} \circ \Psi_{\varepsilon,\Delta \tau}\,
({\bf z}_\varepsilon )-  
\Phi_{\Delta \tau,\tilde{H}_{\Delta \tau}} ({\bf z}_\varepsilon)\|<
d_1\,e^{-d_2/\Delta \tau}
\end{equation}
where ${\bf z}_\varepsilon = \Psi_{\varepsilon,\Delta \tau}^{-1}({\bf z})$ 
denotes the transformed variable, 
$\Phi_{\tau,\tilde H_{\Delta \tau}}$ is the exact flow map
of a Hamiltonian system with transformed Hamiltonian 
\begin{equation}
\tilde{H}_{\Delta \tau} = H_{\Delta \tau} \circ 
\Psi_{\varepsilon,\Delta \tau} = K + \varepsilon {G}_{\Delta \tau} +
e^{-d_3/\varepsilon} {R}_{\Delta \tau},
\end{equation}
$H_{\Delta \tau}$ is the exponentially accurate modified
Hamiltonian for $M_{\Delta t}$ ({\emph i.e.}, (\ref{eq:modhamest}) holds), 
and $d_i>0$, $i=1,2,3$, are appropriate constants.
The transformation is chosen such that
\begin{equation}
\{K,{G}_{\Delta \tau}\} = 0.
\end{equation}
Furthermore, as $\Delta \tau\to0$, 
\begin{equation}
{G}_{\Delta \tau}\to {G}, \qquad
{R}_{\Delta \tau} \to {R},
\end{equation}
where ${G}$ and ${R}$ are given in theorem \ref{normal form theorem}.
\end{theorem}
\begin{proof}
  The proof combines normal form estimates with backward error
  analysis for symplectic methods as developed in \cite{Reich98} 
  for adiabatic invariants in highly oscillatory mechanical
  systems. In particular, backward error analysis implies 
  that the trajectory of the symplectic method $M_{\Delta \tau}$ 
  stays $\Delta \tau$-exponentially close to the exact 
  solution of a system with a
  modified Hamiltonian $H_{\Delta \tau}$ consisting of $H_0$ perturbed by an
  $\mathcal{O}(\varepsilon \Delta \tau^p)$ correction, $p\ge 1$ the order
  of the method. We then apply theorem
  \ref{normal form theorem} to the modified Hamiltonian $H_{\Delta \tau}$.
\end{proof}
Theorem \ref{backward error analysis} implies that 
if $\p^0$ is initialised at zero in the transformed
coordinates, then $\p^n$ will stay exponentially close to zero
in the transformed coordinates over exponentially long times. 
This statement agrees with general results
on the preservation of adiabatic invariants under symplectic time-stepping
methods. See, for example, \cite{Reich98} and \cite{LR2005}.
We will make use of this property of symplectic methods in the following
section. 

Another immediate (more practical) 
conclusion from theorem \ref{backward error analysis} 
is the following 

\begin{corollary} \label{corollary2}
Numerical trajectories computed with a symplectic method
such as (\ref{eq:time_stepping}) satisfy
\begin{equation} \label{eq:numerical-balance}
{\bf p}^n = -\varepsilon J_2 \nabla_{\bf q} V({\bf q}^n) + {\cal
O}(\varepsilon^2), 
\end{equation}
over exponentially many time-steps $n$ provided the conditions of 
theorem \ref{backward error analysis} are satisfied and 
(\ref{eq:numerical-balance}) holds at $n=0$.
\end{corollary}

\begin{proof}
The proof is similar to that of corollary \ref{corollary1}.
\end{proof}

\subsection{A numerical experiment}
\label{modelnumerics}

We consider the particular potential energy function
\begin{equation}
V({\bf q}) = q_x + \frac{1}{2} \exp(-(q_x^2 + q_y^2)), \qquad
{\bf q} = (q_x,q_y)^T,
\end{equation}
and initial conditions
\begin{equation}
q_x(0) = 0, \quad q_y(0) < -10, \quad p_x(0) = 0, \quad
p_y(0) = \varepsilon ,
\end{equation}
{\it i.e.}, $V({\bf q}(0)) \approx 0$. Simulations are run for
long enough such that $q_y(T) > 10$ at final time $T$. The
nonlinear dynamics reduces to a linear system
\begin{equation}
\dot{q}_x = 0, \qquad \dot{q}_y = \varepsilon
\end{equation}
at both $\tau = 0$ and $\tau = T$, {\it i.e.}, to a system
in perfect balance. More precisely, according to
the theoretical results of sections \ref{exponential}
\& \ref{backward}, this statement
should be true at final time up to terms exponentially small in
$\varepsilon$.

\begin{figure}[hbt!]
\begin{center}
\scalebox{0.65}{\includegraphics{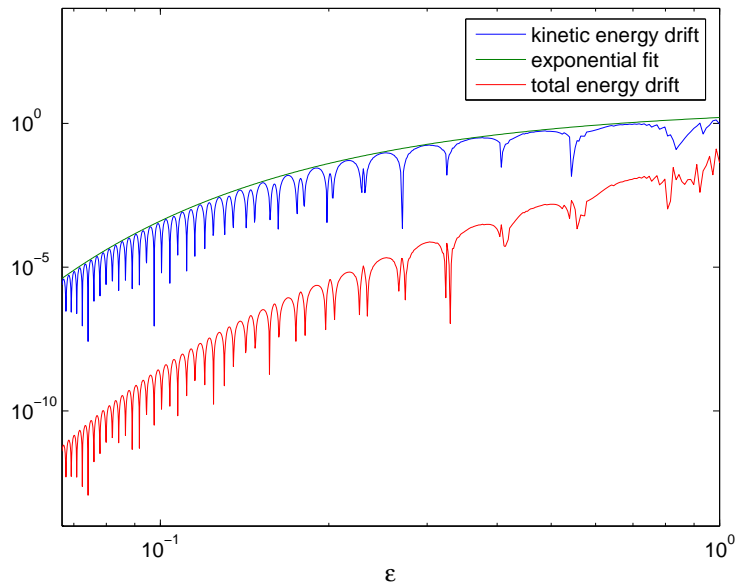}}
\caption{Drift in total energy, kinetic energy, and exponential
fit with $4 \exp(-0.92/\varepsilon)$ as a function of $\varepsilon$.}
\label{figure1}
\end{center}
\end{figure}

We conduct a numerical experiment to verify this
claim. We implement the symplectic method (\ref{eq:time_stepping}) with
step size $\Delta \tau = \varepsilon^2$ and vary $\varepsilon$
in the interval $\tau \in [1/15,1]$. The drift in total energy
$\Delta E = |H({\bf q}(0),{\bf p}(0)) - H({\bf q}(T),{\bf p}(T))|$ and
in kinetic energy $\Delta K = |K({\bf p}(0)) -
K({\bf p}(T))|$ is monitored as a function of $\varepsilon$ and
the results are displayed in figure \ref{figure1}. The drift
in total energy is solely caused by the symplectic time stepping method.
The choice of $\Delta \tau$ ensures that this drift is orders of
magnitude smaller than the computed drift in kinetic energy. In other
words, finite step size effects on the generation of unbalanced motion
can be neglected. The anticipated exponential dependence of the drift
in kinetic energy is nicely confirmed and
\begin{equation}
\Delta K \le 4 \,\exp(-0.92/\varepsilon)
\end{equation}
provides an excellent bound confirming theorem \ref{backward error analysis}
with $d_3 = 0.92$. 

%
\section{Particle methods for the shallow-water equations}
\label{particles}

Throughout section \ref{modelsystem}, 
we have assumed that the fluid depth $\mu$ is a given
function of space and time and derived equations of motion for
a single `fluid parcel'. In this section, we consider
a numerical approximation for $\mu$ of the form
\begin{equation} \label{eq:sph1}
\mu(t,{\bf x}) = \sum_{i=1}^N m_i \,\psi(\|{\bf x} - {\bf q}_i(t)\|)
\end{equation}
in terms of $N$ moving `fluid parcels',
where ${\bf q}_i(t) \in \mathbb{R}^2$ denotes the location of 
the $i$th fluid parcel at time $t$ with mass $m_i$ and shape function
$\psi(r)\ge 0$. The approximation (\ref{eq:sph1}) provides the
typical starting point for a particle method such as smoothed
particle hydrodynamics (SPH), which was first proposed 
in \cite{lucy77,gingold77} for general fluid dynamics and
for the SWEs (\ref{eq:u})-(\ref{eq:mu}) in \cite{salmon83}.

Each `fluid parcel' moves under the Newtonian equations of motion
\begin{eqnarray} \label{eq:sph2}
{\rm Ro} \frac{d}{dt} {\bf p}_i &=& J_2 {\bf p}_i - \frac{\rm B}{\rm Ro}
\nabla_{\bf x} \mu(t,{\bf x})|_{{\bf x} = {\bf q}_i},\\
\frac{d}{dt} {\bf q}_i &=& {\bf p}_i \label{eq:sph3}
\end{eqnarray}
The equations (\ref{eq:sph1})-(\ref{eq:sph3}) form a closed set
of equations and provide an approximation to 
(\ref{eq:nd_u})-(\ref{eq:nd_mu}). (See \cite{frank02sp,frank02bit,FrRe03,bhr} 
for a numerically robust implementation of a particle method for the
SWEs and their geometric properties.)

After setting ${\rm B} = {\rm Ro} = \varepsilon$
and a rescaling of time, equations (\ref{eq:sph2})-(\ref{eq:sph3}) become
\begin{eqnarray} \label{eq:sph2s}
\dot{\bf p}_i &=& J_2 {\bf p}_i - \varepsilon
\nabla_{\bf x} \mu({\bf x},t)|_{{\bf x} = {\bf q}_i},\\
\dot{\bf q}_i &=& {\bf p}_i \label{eq:sph3s}.
\end{eqnarray}
Without restriction of generality, we may also assume that
all `fluid parcels' carry a constant mass $m_i = \delta$.
The equations (\ref{eq:sph1}), (\ref{eq:sph2s})-(\ref{eq:sph3s}) 
are Hamiltonian with Hamiltonian function
\begin{equation} \label{eq:sph_ham}
H_{\rm sph}({\bf z}) = K({\bf p}) + \varepsilon V({\bf q}) = 
\sum_{i=1}^N \frac{1}{2} \|{\bf p}_i\|^2 + \frac{\varepsilon\,\delta}{2} 
\sum_{i,j}^N  \psi(\|{\bf q}_i - {\bf q}_j\|),
\end{equation}
structure matrix
\begin{equation} \label{eq:sph_sm}
J = \left( \begin{array}{cc} J_{2N} & -I_{2N} \\
I_{2N} & 0_{2N} \end{array} \right) \in \mathbb{R}^{2N \times 2N},
\end{equation}
phase space variable ${\bf z} = ({\bf q}^T,{\bf p}^T)^T
\in \mathbb{R}^{4N}$, and 
${\bf q} = ({\bf q}_1^T,\ldots,{\bf q}_N^T)^T$, 
 ${\bf p} = ({\bf p}_1^T,\ldots,{\bf p}_N^T)^T$.
Note that a finite fluid depth $\mu$  implies that the particle 
masses $m_i = \delta$ approach zero as the number of particles
$N \to \infty$. More precisely, $\delta \,N \approx const.$ as
$N\to \infty$.

%
\subsection{Preservation of geostrophic motion}
\label{geostroph}

We first note that the Hamiltonian system defined by
(\ref{eq:sph_ham})-(\ref{eq:sph_sm}) with $N=1$ fits exactly into the
framework of section 2. Furthermore, as first observed in
\cite{cotter04}, the normal form theory developed in section 2 can be
generalised to $N>1$ without much difficulty.  The key idea goes back
to an observation in \cite{ben87:realisation}.  Specifically, it turns
out that the normal form theory for a system with many degrees of
motion with identical fast frequency proceeds much along the lines of
the theory for a single fast degree of motion. Most importantly, the
exponential estimates do not depend on the number of degrees of
freedom and we obtain:

\begin{theorem}
\label{normal form theorem sph}
We consider a fixed number of particles $N$. 
Let the shape function $\psi$ be analytic in some compact subset 
of $\mathbb{C}^{2N}$ containing the solution trajectory. 
Then there exists a near-identity
change of coordinates $\Psi_{\varepsilon}$ such that
\begin{equation} \label{eq:sph_estimate}
\tilde{H} = 
H_{\rm sph}\circ\Psi_{\varepsilon} = K+\varepsilon \,G
+e^{-c/\varepsilon}\,R
\end{equation}
where
\begin{equation}
\{G,K\}=0
\end{equation}
and $c>0$ is some constant.
\end{theorem}
\begin{proof}
We apply theorem \ref{normal form theorem} in slightly modified form
taking into account that $N>1$. We also note that a bound on $\psi$ 
immediately implies a bound on the potential
energy $V$ {\it via} $|V({\bf q})| \le \delta N |\psi |_\infty$. 
\end{proof}

The estimate (\ref{eq:KE_drift}) hence also applies to 
the particle approximation (\ref{eq:sph2s})-(\ref{eq:sph3s}) 
and geostrophic balance is maintained over exponentially
long time intervals. Furthermore, theorem \ref{backward error analysis}
and corollary \ref{corollary2} immediately carry over
to the particle method. Note that (\ref{eq:numerical-balance})
gets replaced by
\begin{equation}
{\bf p}^n_i = -\varepsilon J_2 \nabla_{\bf x} 
\mu({\bf x},t)|_{{\bf x} = {\bf q}_i} +{\cal O}(\varepsilon^2)
\end{equation}
for $i = 1,\ldots,N$. Recall that we consider the limit
$\varepsilon \to 0$ for $N$ fixed.

It is important to keep in mind that theorem \ref{normal form theorem sph}
does {\it not} apply to most practical implementations of particle method since
the basis functions are typically not analytic. In those case the
exponentially small term in (\ref{eq:sph_estimate}) needs to be replaced
by a polynomial expression of the form $(c\varepsilon)^k$, where the
integer $k\ge 2$ and the constant $c>0$ 
depend on the smoothness of the basis function $\psi$. 
This is the case, for example, for the Hamiltonian-Particle-Mesh
(HPM) method (see \cite{frank02sp}), which will be used in the
following section to conduct a number of numerical experiments.

%
\subsection{Numerical experiments}
\label{numerics}

In this section, two numerical examples are given which illustrates
theorems \ref{backward error analysis} and \ref{normal form theorem sph}
for practical implementations of a particle method.

\subsubsection{Exchange of kinetic energy} 

To show that particles may exchange energy as long as
the total kinetic energy is preserved, consider a model consisting of
2 interacting particles moving on a predefined background height
field $\overline{\mu}(\tau,\MM{x})$. In this experiment, we use the Hamiltonian
Particle-Mesh (HPM) method for spatial discretisation (see \cite{frank02sp})
with the total fluid depth at a grid point ${\bf x}_{mn}$ 
given by
\begin{equation}
\mu_{mn}(\tau) = \sum_{i=1}^2\psi_{mn}({\bf q}_i(\tau))
+\overline{\mu}(\tau,{\bf x}_{mn}).
\end{equation}
The basis functions $\psi_{mn}$ are normalized cubic B-splines
centered about grid points ${\bf x}_{mn}$.

We multiplied the potential energy by an analytic function
$g(\tau)$ with $g(\tau)\to 0$ as $\tau \to\pm\infty$ so that the Hamiltonian
takes the form
\begin{equation} 
H_{\rm hpm}({\bf z}) = 
\sum_{i=1}^2 \frac{1}{2} \|{\bf p}_i\|^2 + 
\frac{\varepsilon}{2}g(\tau) 
\sum_{mn}\sum_{i=1}^2\psi_{mn}({\bf q}_i)
\left(\sum_{j=1}^2\psi_{mn}({\bf q}_j)+2\overline{\mu}(\tau,{\bf x}_{mn})
\right).
\end{equation}
The reason for the introduction of the function $g(\tau)$ is that when
$g(\tau)\approx 0$ at the beginning and the end of the simulation, the
Hamiltonian is already in normal form and we can measure the
amount of fast rotational energy exactly (this technique was first
introduced in \cite{g.99:from_hamiltonian} in an experiment to measure
the exchange of oscillatory energy between two
colliding molecules).

\begin{figure}[htp]
\begin{center}
\scalebox{0.75}{\includegraphics{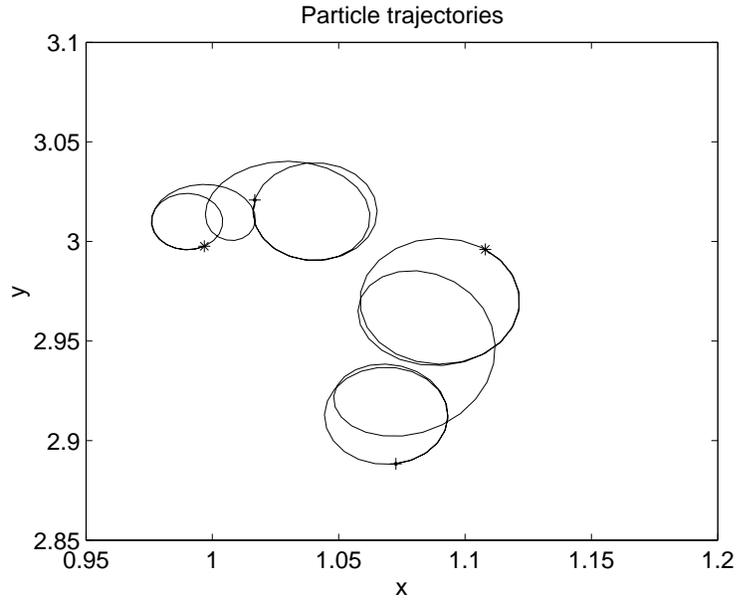}}
\end{center}
\caption{\label{two_part_traj} Plots showing particle trajectories for
  a system of 2 interacting HPM particles moving on a prescribed
  background height field. The initial positions are marked with a `*'
  and the end positions are marked with a `+'. The kinetic energy for
  each particle is shown in figure \ref{two_part_kin}.}
\end{figure}
\begin{figure}[htp]
\begin{center}
\scalebox{0.75}{\includegraphics{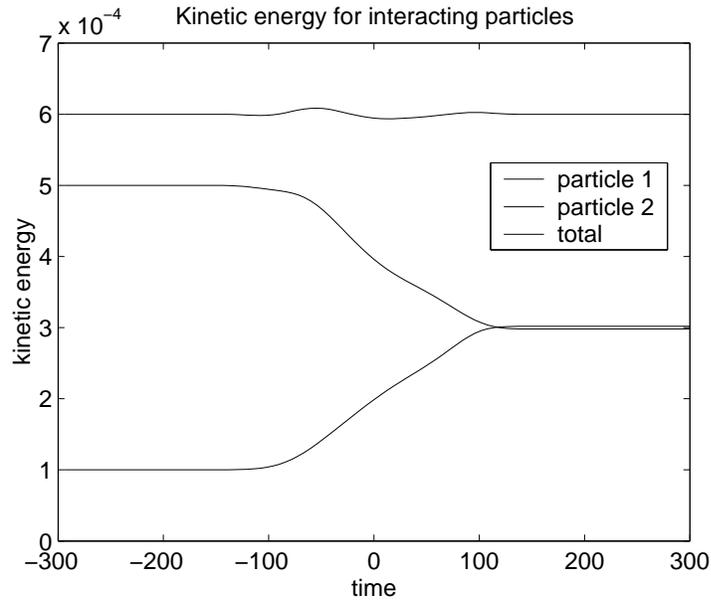}}
\end{center}
\caption{\label{two_part_kin} Plots showing kinetic energies for 2
  interacting HPM particles moving on a prescribed background height
  field. The total kinetic energy returns at the end of the experiment
  to a value very close to the starting value. However the 2 particles
  do exchange kinetic energy.}
\end{figure}
A graph showing kinetic energy against time during the experiment, and
a plot of the trajectories of the two particles is shown in figure
\ref{two_part_traj}.  This example illustrates that, although the
total kinetic energy $K$ is almost invariant, the particles exchange
energy through the interaction potential as shown in figure
\ref{two_part_kin}. This exchange is permitted because the frequencies
of oscillation of the two particles are the same for $\varepsilon =0$. 

\subsubsection{Preservation of geostrophic balance under the
HPM method}

As already eluded to, 
theorem \ref{normal form theorem sph} does not directly apply to
the HPM method \cite{frank02sp} since the basis functions are not
analytic. However, one would nevertheless
expect `good' preservation of geostrophic balance in the semi-geostrophic
scaling limit. To test this hypothesis numerically, we repeated the
shear flow instability simulation of \cite{frank02sp} under slightly
modified initial conditions such that the associated Rossby and Burger
number satisfy ${\rm Ro} \approx {\rm Bu} \approx 0.1$. The HPM
discretized shallow water equations are simulated using $N=262144$ particles
over a domain $(x,y) \in [0,2\pi]^2$, a smoothing length of $\alpha =
0.2015$, and a time step of $\Delta t = 1/36$. Note that one unit of time
corresponds to one day and that the Rossby radius of deformation
corresponds to $L_R \approx 0.5$ in dimensionless variables. The initial
momenta are in geostrophic balance, {\it i.e.},
\begin{equation}
{\bf p}(0) = -\varepsilon J_{2N}\nabla_{\bf q} V({\bf q}).
\end{equation}

\begin{figure}[hbt!]
\begin{center}
\scalebox{0.65}{\includegraphics{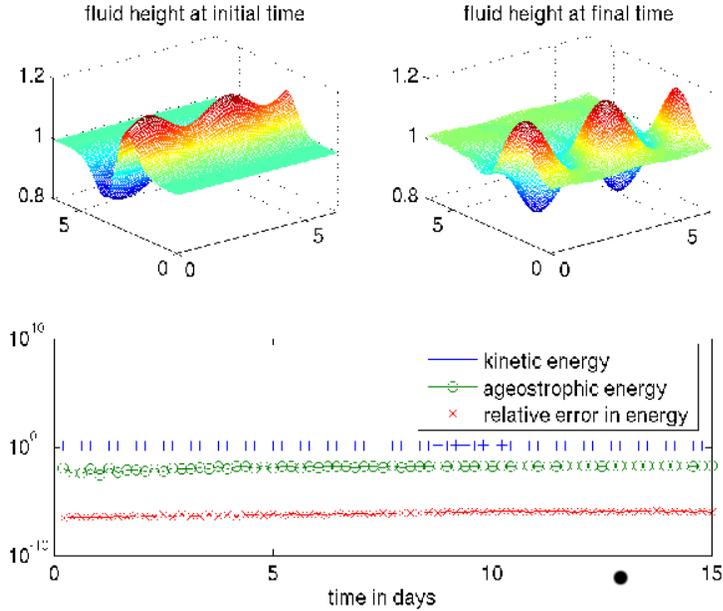}}
\caption{Drift in total energy, kinetic energy, and exponential
fit with $4 \exp(-0.92/\varepsilon)$ as a function of $\varepsilon$.}
\label{figure3}
\end{center}
\end{figure}

Simulation results can be found in figure \ref{figure3}, where
we display the initial and final layer depth ($t=15\,\mbox{days}$).
We also monitor the total kinetic energy $K({\bf p})$, the
total kinetic energy in the ageostrophic momentum components
\begin{equation}
{\bf p}_{\rm ag} = {\bf p} - {\bf p}_{\rm gs}, \qquad
{\bf p}_{\rm gs} = -\varepsilon J_{2N}\nabla_{\bf q} V({\bf q}),
\end{equation}
{\it i.e.}, $K({\bf p}_{\rm ag})$, and the relative error in total
energy $H_{\rm hpm}$. We find that the relative error in energy remains
very small. This is due to the symplectic nature of the time stepping
procedure. We also observe that $K({\bf p}_{\rm ag})$ remains
nearly level and much smaller than the kinetic energy $K({\bf p})$.
This indicates that geostrophic balance is well preserved throughout
the simulation even though $\varepsilon \approx 0.1$ is not very small
and the basis functions $\psi_{mn}$ are not analytic.

%
\section{Summary and outlook}
\label{summary}

We introduced the model system of a single fluid particle moving in a
rotating shallow-water system and gave an exponentially-accurate
normal form theorem valid in the semi-geostrophic scaling limit. This
normal form theorem gives a coordinate change 
for which the kinetic energy stays almost
invariant for very long time intervals. When the kinetic energy in
transformed variables is
initially zero (or close to zero) we showed that it is possible to derive an
exponentially-accurate slow equation which is equivalent to the 
symplectic slow equation obtained by variational asymptotics in the 
semi-geostrophic limit.

We showed that this result may be extended to numerical schemes for
solving the SWEs based on particle methods by
making use of backward error analysis, and then illustrated this with
a numerical experiment designed to ``capture'' the exact amount of 
energy exchanged between the slow and fast dynamics.

The extension of this work to the full rotating shallow-water PDE
(\ref{eq:u})-(\ref{eq:mu})
requires bounds to be obtained on the gradients of the height
field. As the SWEs are known to develop shocks it
seems unlikely that an exponentially-accurate result will be valid for
long times unless it is possible to obtain some extra regularity when
the fast dynamics has very little energy. It may be possible to obtain
estimates for regularised equations where a smoothing operator is
applied to the height gradient in the momentum equation, or for the
family of $\alpha$-regularised equations such as
shallow-water-$\alpha$ \cite{holm99}.

Finally we suggest one possibly fruitful application of this work in
the calculation of the dispersion of passive tracers using Lagrangian
particles. Typically these calculations are performed by interpolating
a given (gridded) velocity field ${\bf u}_{mn}(t)$ from the grid to the 
whole domain, and solving the
system of ODEs
\begin{equation}
\frac{d {\bf q}_i}{dt} 
= \sum_{mn} \psi_{mn}({\bf q}_i)\,{\bf u}_{mn}(t), 
\qquad i=1,\ldots,N.
\end{equation}
This method can often experience problems with artificial clumping of
particles when the trajectories are integrated over long time
intervals. An alternative approach would be to solve equations
(\ref{eq:ddot_x})-(\ref{eq:ddot_y}) with a given fluid depth (or
geopotential) $\mu$. The results of this paper show that if the system is
integrated with a suitable symplectic method then the tracer particles
will keep balanced trajectories in the semi-geostrophic limit.

%
\section*{Acknowledgement}

We thank Georg Gottwald and Marcel Oliver for stimulating discussions
on the relation of the Hamiltonian normal form theory, as proposed in
\cite{cotter04}, and the asymptotic variational approach of
\cite{oliver05,gottwald1,gottwald2}.

%
\bibliography{survey}
%


%
\section*{Appendix. Calculation of slow equation in symplectic form
  to second-order}
%

In this appendix we explicitly calculate the change of coordinates and
the slow equation to $\mathcal{O}(\epsilon^2)$ using the iterative
procedure given in section \ref{exponential}.

The Hamiltonian vector field given by the kinetic energy 
$K = \|{\bf p}\|^2/2$ is 
\[
\dot{\p} =  J_2\p, \qquad \dot{\q} = \p,
\]
and the associated time-$\tau$ flow map $\Phi_{\tau,K}$ is
\begin{eqnarray*}
\p & \mapsto & e^{J_2\tau}\p, \\
\q & \mapsto & J_2(1-e^{J_2\tau})\p + \q.
\end{eqnarray*}

First we need to calculate the transformed Hamiltonian after
the first-order transformation.
\begin{eqnarray*}
H_1 &=& H_0\circ\Phi_{1,\varepsilon F_1}, \\
&=& H_0 + \varepsilon\{H_0,F_1\} + \frac{\varepsilon^2}{2}\{\{H_0,F_1\},F_1\}
+ \mathcal{O}(\varepsilon^3), \\
&=& K + \varepsilon V + \varepsilon\{K,F_1\} + \varepsilon^2\{V,g_1\}
+\frac{\varepsilon^2}{2}\{\{K,F_1\},F_1\} + \mathcal{O}(\varepsilon^3), \\
&=& K + \varepsilon\bar{V} + \frac{\varepsilon^2}{2}\{\bar{V}+V,F_1\} +
\mathcal{O}(\varepsilon^3), \\
\end{eqnarray*}
where we have used $\{K,g_1\} = \bar{V}-V$, $\{K,\bar{V}\} = 0$, and
\begin{eqnarray*}
\bar{V}({\bf z}) & = & \frac{1}{T}\int_0^T d\tau \,V\circ
\Phi_{\tau,K}({\bf z}), \\
& = & \frac{1}{T}\int_0^T d\tau\, V\left(J_2 \left(1-e^{J_2\tau}
\right)\p+\q\right), 
\end{eqnarray*}
as well as
\begin{eqnarray*}
F_1({\bf z}) & = & \frac{1}{T}\int_0^T d\tau\,\tau 
(V-\bar{V})\circ\Phi_{\tau,K}({\bf z}), \\
& = & \frac{1}{T}\int_0^T d\tau \, \tau \left(
V\left(J_2 \left(1-e^{J_2\tau} \right)\p+\q\right)-\int_0^T d\tau'\, 
V\left(J_2 \left(1-e^{J_2(\tau+\tau')} \right)\p+\q\right)
\right).
\end{eqnarray*}
Recall that ${\bf z} = ({\bf p}^T,{\bf q}^T)^T$ and $T = 2\pi$. 

The first-order slow equation is given by:
\[
J_2\dot{\q} = \varepsilon\nabla_{\q}\bar{V}\left({\bf p},{\bf q}\right)
|_{\p={\bf 0}} 
= \varepsilon\nabla_{\bf q} V(\q).
\]
To calculate the slow equation at the next order, we write the
transformed Hamiltonian after the first-order transformation as
\begin{eqnarray*}
H_2 &=& K + \varepsilon\bar{V} + \varepsilon^2\bar{R_1} + \mathcal{O}(\varepsilon^3),
\end{eqnarray*}
where
\[
R_1 = \frac{1}{2}\{\bar{V}+V,F_1\}.
\]
Note that for a general function $f(\p,\q)$, 
\[
\nabla_{\q}\bar{f}(\p,\q)|_{\p={\bf 0}} =
\frac{1}{T}\int_0^T d\tau\, 
\nabla_{\q} \left\{f \left(e^{J_2\tau}\p,J_2 \left(1-e^{J_2\tau}\right)\p
+\q\right)\right\}
|_{\p={\bf 0}} = \nabla_{\q} f({\bf 0},\q),
\]
so all that remains is to calculate $R_1$ and the gradient with respect 
to $\q$ for $\p= {\bf 0}$. 

We have
\[
\{\bar{V},F_1\} = \nabla_{\p}\bar{V}\cdot J_2\nabla_{\p}F_1
+ \nabla_{\q}\bar{V}\cdot\nabla_{\p}F_1 -
\nabla_{\p} \bar{V}\cdot \nabla_{\q}F_1,
\]
and
\[
\{V,F_1\} = \nabla_{\q}V\cdot\nabla_{\p}F_1.
\]
Now we need to calculate these various gradients:
\begin{eqnarray*}
\nabla_{\q}F_1({\bf z})|_{\p={\bf 0}} & = & \frac{1}{T}\int_0^T
d\tau\, \tau \nabla_{\q} (V-\bar{V})\circ\Phi_{\tau,K}({\bf z})
|_{\p={\bf 0}} , \\
& = & \nabla_{\q}(V({\bf q})-V({\bf q})) = {\bf 0},
\end{eqnarray*}
and
\begin{eqnarray*}
\nabla_{\p}F_1({\bf z})|_{\p={\bf 0}} & = & \frac{1}{T}\int_0^T
d\tau\, \tau \Bigg[
J_2 \left(1-e^{J_2\tau}\right)\nabla_{\bf q} V \left(J_2 \left(1-e^{J_2\tau}
\right)\p+\q \right)|_{\p={\bf 0}} - \\
& & \qquad\qquad \int_0^T d\tau' \,J_2 \left(1-e^{J_2(d\tau+d\tau')}\right)
\nabla_{\bf q}V \left(J_2\left(1-e^{J_2\tau}\right)\p+\q\right)|_{\p={\bf 0}}
\Bigg], \\
& = & \left(
\frac{-1}{T}\int_0^T d\tau \, \tau \left[
J_2 \left(1-e^{J_2\tau} \right)-
\int_0^T d\tau' \, J_2\left(1-e^{J_2(d\tau+d\tau')}\right) 
 \right]
\right)\nabla_{\bf q} V(\q), \\
& = & \left(
\frac{1}{T}\int_0^T d\tau\,
\tau \left[J_2 \left(1-e^{J_2\tau}\right)-J_2\right]
\right)\nabla_{\bf q} V(\q), \\
& = & \frac{-1}{T}\int_0^T 
d\tau\,J_2\tau e^{J_2\tau}  \nabla_{\bf q} V(\q), \\
& = & \frac{-1}{T}\int_0^T d\tau 
\, \tau\frac{d}{d\tau}e^{J_2\tau}  \nabla_{\bf q} V(\q), \\
& = & -\nabla_{\bf q} V(\q).
\end{eqnarray*}
We also obtain
\begin{eqnarray*}
\nabla_{\q}\bar{V}({\bf z})|_{\p={\bf 0}} & = & \frac{1}{T}\int_0^T
d\tau\, \nabla_{\bf q} V
\left(J_2 \left(1-e^{J_2\tau}\right)\p+\q\right))|_{\p={\bf 0}}, 
\\& = & \nabla_{\bf q} V(\q).
\end{eqnarray*}
and
\begin{eqnarray*}
\nabla_{\p}\bar{V}({\bf z})|_{\p={\bf 0}} & = & \frac{-1}{T}\int_0^T
d\tau\,\left(1-e^{-J_2\tau}\right)J_2\nabla_{\q}V\left(J_2\left(1
-e^{J_2\tau}\right)\p+\q\right)|_{\p=0}, \\
& = & -J_2\nabla_{\bf q} V(\q).
\end{eqnarray*}
Hence we get
\[
\{\bar{V},F_1\}|_{\p={\bf 0}} = -\left(J_2\nabla_{\bf q} 
V\right)\cdot \left(J_2\nabla_{\bf q} V\right) 
- \nabla_{\bf q} V \cdot \nabla_{\bf q} V
= 0,
\]
and
\[
\{V,F_1\}|_{\p={\bf 0}} = 
-\nabla_{\bf q} V\cdot\nabla_{\bf q} V = -\|\nabla_{\bf q} V\|^2,
\]
and so $R_1 = -\frac{1}{2}\|\nabla_{\bf q} V\|^2$.
The slow equation in transformed coordinates is then
\[
J_2\dot{\bf q}_\varepsilon = \varepsilon\nabla_{\bf q} V({\bf q}_\varepsilon) 
- \frac{\varepsilon^2}{2}\nabla_{\bf q} \|\nabla_{\bf q} 
V({\bf q}_\varepsilon)\|^2 + \mathcal{O}(\varepsilon^3).
\]
We note that this is the same slow canonical equation obtained using 
variational asymptotics in \cite{gottwald1,gottwald2} and, in particular,
yields the `large-scale semi-geostrophic' equation (\ref{eq:igsa}).

\end{document}